\documentclass[12pt]{article}
\usepackage{amssymb,amsmath}
\def\D{\Delta}
\def\UU{{\cal U}}
\def\WW{{\cal L}}
\def\LL{{\cal L}}
\def\Der{{\rm Der}}
\def\VV{\mathcal {V}}
\def\Inn{{\rm Inn}}
\def\Ker{{\rm Ker}}
\def\Im{{\rm Im}}
\def\v{\varphi}
\def\ssc{\scriptscriptstyle}
\def\cl{\centerline}
\def\ma{\mathbb}
\def\rar{\rightarrow}
\def\bs{\backslash}
\def\vs{\vspace*}
\def\even{{\bar0}}
\def\odd{{\bar1}}
\def\bo{{\bf1}}
\def\ni{\noindent}
\def\N{\mathbb{N}{\ssc\,}}
\def\Z{\mathbb{Z}{\ssc\,}}
\def\C{\mathbb{C}{\ssc\,}}

\def\adddot{$\!\!\!${\bf.}\ \ }
\def\QED{\hfill$\Box$}

\numberwithin{equation}{section}
\newtheorem{theo}{Theorem}[section]
\newtheorem{defi}[theo]{Definition}

\newtheorem{lemm}[theo]{Lemma}
\newtheorem{prop}[theo]{Proposition}
\newtheorem{clai}{Claim}

\newtheorem{coro}[theo]{Corollary}

\begin{document}\baselineskip 18pt
\cl{\large \bf Lie superbialgebra structures on the centerless}
\cl{{\large \bf twisted $N=2$ superconformal algebra}\footnote
{Supported by
NSF grants 10671027, 10825101 of China\\[2pt] \indent\ \,$^{\dag}$Corresponding E-mail:
sd\_junbo@163.com}} \vs{6pt}

\cl{Huanxia Fa, Junbo Li \!$^{\dag}$} \cl{\small Department of
Mathematics, Changshu Institute of Technology, Changshu 215500,
China} \vs{6pt}

{\small
\parskip .005 truein
\baselineskip 3pt \lineskip 3pt

\noindent{{\bf Abstract.} In this paper, Lie superbialgebra
structures on the centerless twisted $N=2$ superconformal algebra
$\LL$ are considered which are proved to be coboundary triangular.
\vs{5pt}

}

}

\vs{6pt}
 \cl{\bf\S1. \
Preliminaries}\setcounter{section}{1}\setcounter{equation}{0}

The superconformal algebras were constructed by Kac (see \cite{K})
and by Ademollo et al. [10] originally but independently, which are
closely related the conformal field theory and the string theory,
and act important roles in both mathematics and physics supplying
the underlying symmetries of string theory. It is well-known that
the $N=2$ superconformal algebras fall into four sectors: the
Neveu-Schwarz sector, the Ramond sector, the topological sector and
the twisted sector. A series of good results have been obtained on
these algebras (see \cite{CK,DB,EG,KL,SS,YS2} and the reference
cited therein). All sectors are closely related to the Virasoro
algebra and the super-Virasoro algebra which play great roles in the
two-dimensional conformal filed.

The notion of Lie bialgebras was introduced in 1983 by Drinfeld (see
\cite{D1,D2}) during the process of investigating quantum groups.
Then there appeared several papers on Lie bialgebras and Lie
superbialgebras (e.g., \cite{M,N,NT,YS1,YS2}). In
\cite{M}--\cite{NT}, the Lie bialgebra structures on Witt and
Virasoro algebras were investigated, which are shown to be
triangular coboundary. Moreover, the Lie bialgebra structures on the
one-sided Witt algebra were completely classified. In
\cite{YS1,YS2}, the Lie superbialgebra structures on the generalized
super-Virasoro algebra and Ramond $N=2$ superconformal algebra were
investigated. In this paper, we shall study the Lie super-bialgebra
structures on the centerless twisted $N=2$ superconformal algebra,
which is proved to be coboundary triangular.

Firstly, let us recall some related definitions. Let
$\LL=\LL_\even\oplus \LL_\odd$ be a vector space over the complex
number field $\C$. If $x\in\LL_{[x]}$, then we say that $x$ is
homogeneous of degree $[x]$ and we write ${\rm deg}x=[x]$. Denote by
$\tau$ the {\it super-twist map} of $\LL \otimes \LL$, i.e.,
\begin{eqnarray*}
\tau(x\otimes y)= (-1)^{[x][y]}y\otimes x,\ \ \,\forall\,\,x,y\in
\LL.
\end{eqnarray*}
For any $n\in\N$, denote by $\LL^{\otimes n}$ the tensor product of
$n$ copies of $\LL$ and $\xi$ the {\it super-cyclic map} cyclically
permuting the coordinates of $\LL^{\otimes3}$, i.e.,
\begin{eqnarray*}
&&\xi=({\bf1}\otimes\tau)\cdot(\tau\otimes\bo): \,x_{1} \otimes
x_{2} \otimes x_{3}\mapsto (-1)^{[x_1]([x_2]+[x_3])}x_{2} \otimes
x_{3} \otimes x_{1},\ \ \forall\,\,x_i\in \LL,\,i=1,2,3,
\end{eqnarray*}
where $\bo$ is the identity map of $\LL$. Then the definition of a
Lie superalgebra can be described in the following way: A {\it Lie
superalgebra} is a pair $(\LL,\v)$ consisting of a vector space
$\LL=\LL_\even\oplus \LL_\odd$ and a bilinear map $\v :\LL
\otimes\LL \to\LL$ satisfying:
\begin{eqnarray}\label{cond1}\begin{array}{lll}
&&\v(\LL_{\bar{i}},\LL_{\bar{j}})\subset \LL_{\bar{i}+\bar{j}},\vs{4pt}\\
&&\Ker(\bo \otimes\bo -\tau) \subset \Ker\,\v,\vs{4pt}\\
&&\v \cdot (\bo  \otimes \v ) \cdot (\bo \otimes\bo \otimes\bo  +
\xi +\xi^{2})=0.
\end{array}\end{eqnarray}
Meanwhile, the definition of a Lie super-coalgebra can be described
in the following way: A {\it Lie super-coalgebra} is a pair
$(\LL,\D)$ consisting of a vector space $\LL=\LL_\even\oplus
\LL_\odd$ and a linear map $\D: \LL \rar \LL \otimes \LL$
satisfying:
\begin{eqnarray}\label{cond2}\begin{array}{lll}
&&\D(\LL_{\bar{i}})\subset
\mbox{$\sum\limits_{\bar{j}\in{\Z_2}}$}\LL_{\bar{j}}\otimes
\LL_{\bar{i}-\bar{j}},\vs{4pt}\\
&& \Im\,\D \subset \Im(\bo \otimes\bo - \tau),\vs{4pt}\\
&&(\bo \otimes\bo \otimes\bo  + \xi +\xi^{2}) \cdot (\bo  \otimes
\D) \cdot \D =0.
\end{array}\end{eqnarray}
Now one can give the definition of a Lie super-bialgebra, which is
a triple $(\LL, \v, \D )$ satisfying:
\begin{eqnarray*}
&&\mbox{(i)}\ \ (\LL, \v){\mbox{ is a Lie superalgebra}},\\[2pt]
&&\mbox{(ii)}\ \ (\LL, \D){\mbox{ is a Lie super-coalgebra}},\\[2pt]
&&\mbox{(iii)}\ \ \D  \v (x\otimes y)=x \ast  \D y - (-1)^{[x][y]}y
\ast  \D x\ \ \forall\,\,x, y \in \LL,
\end{eqnarray*} where
the symbol ``$\ast$'' means the {\it adjoint diagonal action}
\begin{eqnarray}\label{e-diag}
x\ast  (\mbox{$\sum\limits_{i}$}{a_{i} \otimes b_{i}}) =
\mbox{$\sum\limits_{i}$} ( {[x, a_{i}] \otimes b_{i} +
(-1)^{[x][a_i]}a_{i} \otimes [x, b_{i}]}),\ \ \forall\,\,x, a_{i},
b_{i} \in \LL,
\end{eqnarray}
and in general $[x,y]=\v(x\otimes y)$ for $x,y \in \LL$.

Denote by $\UU(\LL )$ the {\it universal enveloping algebra} of
$\LL$ and $A\bs B=\{x\,|\,x\in A,x\notin B\}$ for any two sets $A$
and $B$. If $r =\mbox{$\sum\limits_{i}$}{a_{i} \otimes b_{i}} \in
\LL \otimes \LL $, then the following elements are in
$\UU(\LL)\otimes \UU(\LL)\otimes \UU(\LL)$
\begin{eqnarray*}
&&r^{12} =\mbox{$\sum\limits_{i}$}{a_{i} \otimes b_{i} \otimes \bo
}=r\otimes\bo,\ r^{23} =\mbox{$\sum\limits_{i}$}{\bo  \otimes a_{i}
\otimes
b_{i}}=\bo \otimes r,\\
&&r^{13}=\mbox{$\sum\limits_{i}$} {a_{i} \otimes \bo  \otimes b_{i}}
=(\bo \otimes\tau)(r\otimes \bo )=(\tau\otimes\bo )(\bo \otimes r),
\end{eqnarray*}
while the following elements are in $\LL\otimes\LL\otimes\LL$
\begin{eqnarray*}
&&[r^{12} ,r^{23}]=\mbox{$\sum\limits_{i,j}$}a_{i} \otimes[
b_{i},a_{j}] \otimes
b_{j},\\
&&[r^{12} ,r^{13}]=\mbox{$\sum\limits_{i,
j}$}(-1)^{[a_j][b_i]}[a_{i}, a_{j}] \otimes
b_{i} \otimes b_{j},\\
&&[r^{13} ,r^{23}]=\mbox{$\sum\limits_{i,j} (-1)^{[a_j][b_i]}$}a_{i}
\otimes a_{j} \otimes[ b_{i}, b_{j}].
\end{eqnarray*}

\begin{defi}\adddot (i) A {\it coboundary super-bialgebra} is a quadruple $(\LL , \v,
\D,r),$ where $(\LL , \v, \D)$ is a Lie super-bialgebra and $r \in
\Im(\bo \otimes\bo  - \tau) \subset \LL \otimes \LL $ such that
$\D=\D_r$ is a {\it coboundary of $r$}, i.e.,
\begin{eqnarray}\label{e-D-r}
\D_r(x)=(-1)^{[r][x]}x\ast r,\ \ \forall\,\,x \in \LL.
\end{eqnarray}
(ii)\ \ A coboundary Lie super-bialgebra $(\LL , \v,\D, r)$ is
called {\it triangular} if it satisfies the following {\it classical
Yang-Baxter Equation}
\begin{eqnarray}\label{e-CYBE}
c(r):=[r^{12}, r^{13}] +[r^{12} , r^{23}] +[r^{13} , r^{23}]=0.\ \ \
{\rm (CYBE)}
\end{eqnarray}
\end{defi}

Let $V=V_{\bar0}\oplus V_{\bar1}$ be an $\LL$-module where
$\LL=\LL_{\bar0}\oplus\LL_{\bar1}$. A $\Z_2$-homogenous linear map
$d:\LL\to V$ is called a {\it homogenous derivation of degree
$[d]\in\Z_2$}, if $d(\LL_i)\subset V_{i+[d]}\
\,(\forall\,\,i\in\Z_2)$,
\begin{eqnarray}\label{e-der}
&&d([x,y])=(-1)^{[d][x]}x\ast d(y)-(-1)^{[y]([d]+[x])}y\ast  d(x),\
\ \forall\,\,x,y\in \LL.
\end{eqnarray}
Denote by $\Der_{\bar{i}}(\LL,V)\ \,(\,i=0,1)$ the set of all
homogenous derivations of degree $\bar{i}$. Then the set of all
derivations from $\LL$ to $V$
$\Der(\LL,V)=\Der_{\bar0}(\LL,V)\oplus\Der_{\bar1}(\LL,V)$. Denote
by $\Inn_{\bar{i}}(\LL,V)\ \,(\,i=0,1)$ the set of {\it homogenous
inner derivations of degree $\bar{i}$}, consisting of $a_{\rm inn},$
$a\in V_{\bar{i}}$, defined by
\begin{equation}
\label{e-inn} a_{\rm inn}:x\mapsto (-1)^{[a][x]}x\ast a,\ \
\forall\,\,x\in \LL.
\end{equation}
Then the set of inner derivations
$\Inn(\LL,V)=\Inn_{\bar0}(\LL,V)\oplus\Inn_{\bar1}(\LL,V)$.

Denote by $H^1(\LL,V)$ the {\it first cohomology group} of $\LL$
with coefficients in $V$. Then
\begin{equation*}
H^1(\LL,V)\cong\Der(\LL,V)/\Inn(\LL,V).
\end{equation*}
An element $r$ in a superalgebra $\LL$ is said to satisfy the {\it
modified Yang-Baxter equation} if
\begin{equation}\label{e-MYBE}
x\ast  c(r)=0,\ \ \forall\,\,x\in \LL.\ \ \ {\rm (MYBE)}
\end{equation}

The {\it centerless twisted $N=2$ superconformal algebra} $\LL$
consists of the Virasoro algebra generators $L_m,\,m\in\Z$,
corresponding to the stress-energy tensor, a Heisenberg algebra
$T_r$, with half-integral $r\in\frac{1}{2}+\Z$, corresponding to the
$U(1)$ current, and the fermionic generators
$G_p,\,p\in\frac{1}{2}\Z$, which are the modes of the two
spin-$\frac{3}{2}$ fermionic fields with the following commutation
relations (see, e.g., \cite{DB}),
\begin{eqnarray}\label{LieB}\begin{array}{lllll}
&& [L_m,L_n]=(m-n)L_{n+m},\vs{6pt}\\
&&[L_m,T_r]=-rT_{r+m},\ \ \ \ \ \ \ \ \ \ \ \ \ \ \ \ [T_r, T_s]=0,\vs{6pt}\\
&&
[L_m,G_p]=(\frac{m}{2}-p)G_{p+m},\ \ \ \ \ \ \ \ [T_r,G_p]=G_{p+r},\vs{6pt}\\
&&[G_p,G_q]=\left\{\begin{array}{ll}(-1)^{2p}2L_{p+q}&\mbox{if \ }p+q\in\Z,\\[6pt]
(-1)^{2p+1}(p-q)T_{p+q}&\mbox{if \
}p+q\in\frac12+\Z,\end{array}\right.
\end{array}\end{eqnarray}
for $m,n\in\Z,\,r,s\in\frac{1}{2}+\Z,\,p,q\in\frac{1}{2}\Z$.
Obviously, $\LL$ is ${\mathbb{Z}}_2$-graded: ${\cal L}={\cal
L}_{\overline{0}}\oplus{\cal L}_{\overline{1}},$ with
\begin{eqnarray}\label{0426a002}
{\cal L}_{\overline{0}}=\mbox{span}_{\ma C}\{
L_m,\,T_r\,|\,m,\,r\in\frac{1}{2}+\Z\},\ \ \ \ {\cal
L}_{\overline{1}}=\mbox{span}_{\ma C}\{G_p\,|\,p\in\frac{1}{2}\Z\}.
\end{eqnarray}
The Cartan subalgebra of ${\cal L}$ is ${\cal H}={\mathbb C}L_0$ and
$\mathcal {W}=$ span$_{\mathbb C} \{ L_m\,|\,m\in{\mathbb Z}\}$ is
the well-known centerless Virasoro algebra.\vs{6pt}

The main result of this paper can be formulated as follows.
\begin{theo}\adddot\label{mainthe}
Every Lie super-bialgebra structure on the centerless twisted $N=2$
superconformal algebra $\LL$ defined in $(\ref{LieB})$ is triangular
coboundary.
\end{theo}

\vs{12pt}

\cl{\bf\S2. \ Proof of the main result}\setcounter{section}{2}
\setcounter{theo}{0} \setcounter{equation}{0} \vs{6pt} The following
result for the non-super case can be found in \cite{NT} while its
super case can be found in \cite{YS2}.
\begin{lemm}\adddot\label{0502m01}
Let $\LL$ be a Lie superalgebra,
$r\in\Im(\bo\otimes\bo-\tau)\subset\WW\otimes\WW$ with $[r]=\bar0$.
Then
\begin{equation}\label{e-p2.1}
(\bo + \xi + \xi^{2}) \cdot (\bo \otimes \D_r) \cdot \D_r (x) = x
\ast c (r),\ \ \forall\,\,x \in \LL.
\end{equation}
Thus $(\LL,[\cdot,\cdot],\D_r)$ is a Lie super-bialgebra if and only
if $r$ satisfies $(MYBE)$ $({\rm see} (\ref{e-MYBE}))$.
\end{lemm}
The following lemma can be obtained by using the similar techniques
of \cite[Lemma 2.2]{SSu}.
\begin{lemm}\adddot \label{lemma2}
Regarding $\LL^{\otimes n}$ as an $\LL$-module under the adjoint
diagonal action of $\LL $, if $r\in\LL^{\otimes n}$ such that $x\ast
r=0,\ \forall\,\,x\in\LL$, then one has $r=0$.
\end{lemm}

As a conclusion of Lemma \ref{lemma2}, one immediately obtains
\begin{coro}\adddot\label{colo}
An element $r \in \Im(\bo\otimes\bo - \tau) \subset \WW \otimes \WW
$ satisfies CYBE in $(\ref{e-CYBE})$ if and only if it satisfies
MYBE in $(\ref{e-MYBE})$.
\end{coro}

\begin{prop}\adddot\label{lemma3}
$\Der(\LL,\VV)=\Inn(\LL,\VV),$ where $\VV=\LL\otimes \LL$,
equivalently, $H^1(\LL,V)=0$.
\end{prop}
\ni{\it Proof.~}\def\b{j}\def\g{k} Note that
$\VV=\oplus_{i\in\frac{1}{2}\Z}\VV_i$ is also $\frac{1}{2}\Z$-graded
with $\VV_i=\sum_{j+k=i} \LL_j\otimes\LL_k$, where
$i,j,k\in\frac{1}{2}\Z$. We say a derivation $d\in\Der(\LL ,\VV)$ is
{\it homogeneous of degree $i\in\frac{1}{2}\Z$} if $d(\VV_j) \subset
\VV_{i +j}$ for all $j\in\frac{1}{2}\Z$. Set $\Der(\LL , \VV)_i
=\{d\in \Der(\LL , \VV) \,|\,{\rm deg\,}d =i\}$ for
$i\in\frac{1}{2}\Z$.

For any $d\in\Der(\LL,\VV)$, $i\in\frac{1}{2}\Z$, $u\in\LL _j$ with
$j\in\frac{1}{2}\Z$, we can write $d(u)=\sum_{k\in\Z}v_k\in \VV$
with $v_k\in \VV_k$, then we set $d_i(u)=v_{i+j}$. Then
$d_i\in\Der(\LL,\VV)_i$ and
\begin{eqnarray}\label{summable}
\mbox{$d=\sum\limits_{i\in\frac{1}{2}\Z} d_i\
\,\mbox{\ where\ }d_i \in \Der(\LL, \VV)_i,$}
\end{eqnarray}
which holds in the sense that for every $u \in\LL $ only finitely
many $d_i(u)\neq 0,$ and $d(u) = \sum_{i \in\Z} d_i(u)$ (we call
such a sum in (\ref{summable}) {\it summable}).

\begin{clai}\adddot\label{clai1}
\rm If \,$i\in\frac{1}{2}\Z\bs\{0\}$, then $d_i\in\Inn(\LL ,\VV)$.
\end{clai}
\ \indent Denote $u=-\frac1i{\ssc\,} d_{i}(L_0)\in \VV_{i}.$ For any
$x_{j}\in \LL_{j},j \in\frac{1}{2}\Z,$  applying $d_{i}$ to
$[L_0,x_{j}]=-j x_{j},$ using $d_{i}(x_j)\in \VV_{i+j}$ and the
action of $L_0$ on $\VV_{i+j}$ is the scalar
$L_0|_{\VV_{i+j}}=-(i+j)$, one has
\begin{equation}\label{equa-add-1}
-(i+j)d_{i}(x_{j}) - (-1)^{[d_{j}][x_{j}]}x_{j}\cdot d_{i}(L_0)=-j
d_{i}(x_{j}),
\end{equation}
i.e., $d_{i}(x_{j})=u_{\rm inn}(x_{j})$, which implies $d_{i}$ is
inner.
\begin{clai}\adddot\label{sub1}
\rm $d_0(L_0)=0.$
\end{clai}
\par

Using (\ref{equa-add-1}) with $i=0$, we obtain $x\ast d_0(L_0)=0,\,\
\forall\,\,x\in\LL_{\b},\,j\in\frac{1}{2}\Z$, which together with
Lemma \ref{lemma2} gives $d_0(L_0)=0$.
\begin{clai}\adddot\label{sub2}
\rm For any $d_0\in\Der_{\bar0}(\WW,\VV)$, replacing $d_0$ by
$d_0-u_{\rm inn}$ for some $u\in \VV_0$, one can suppose
$d_0(\LL)=0.$
\end{clai}

For any $n\in\Z^*$, one can write $d_0(L_{n})$ as
\begin{eqnarray}\label{0502n01}
&&\mbox{$\sum\limits_{i\in\Z}$}a_{n,i}L_{i+ n}\otimes
L_{-i}+\mbox{$\sum\limits_{j\in\Z}$}(d_{n,j}L_{j+ n}\otimes
G_{-j}+e_{n,j}G_{j+n}\otimes
L_{-j})+\mbox{$\sum\limits_{p\in\frac{1}{2}\Z}$}b_{n,p}G_{p+
n}\otimes
G_{-p}\nonumber\\
&&+\mbox{$\sum\limits_{r\in\frac{1}{2}+\Z}$}c_{n,r}T_{r+
n}\otimes
T_{-r}+\mbox{$\sum\limits_{r\in\frac{1}{2}+\Z}$}(f_{n,r}G_{r+
n}\otimes T_{-r}+g_{n,r}T_{r+ n}\otimes G_{-r}),
\end{eqnarray}
for some $a_{n,i},b_{n,p},c_{n,r},d_{n,j},e_{n,j},f_{n,r},g_{n,r}\in
\C$, where the sums are all finite. Noticing that for any
$i,j\in\Z,\,r\in\frac{1}{2}+\Z,\,p\in\frac{1}{2}\Z$, one has
\begin{eqnarray*}
&&L_1\ast(T_{r}\otimes T_{-r})=rT_{r}\otimes T_{1-r}-rT_{1+r}\otimes
T_{-r},\\
&&L_1\ast(T_{r}\otimes G_{-r})=(0.5+r)T_{r}\otimes
G_{1-r}-rT_{1+r}\otimes G_{-r},\\
&&L_1\ast(G_{r}\otimes T_{-r})=rG_{r}\otimes
T_{1-r}-(r-0.5)G_{1+r}\otimes
T_{-r},\\
&&L_1\ast(L_{i}\otimes L_{-i})=(i+1)L_{i}\otimes
L_{1-i}-(i-1)L_{1+i}\otimes
L_{-i},
\end{eqnarray*}
and
\begin{eqnarray*}
&&L_1\ast(L_{j}\otimes G_{-j})=(j+0.5)L_{j}\otimes
G_{1-j}-(j-1)L_{1+j}\otimes G_{-j},\\
&&L_1\ast(G_{j}\otimes L_{-j})=(j+1)G_{j}\otimes
L_{1-j}-(j-0.5)G_{1+j}\otimes L_{-j},\\
&&L_1\ast(G_{p}\otimes G_{-p})=(p+0.5)G_{p}\otimes
G_{1-p}-(p-0.5)G_{1+p}\otimes G_{-p}\,.
\end{eqnarray*}
Denote
\begin{eqnarray*}
M_{1,1}=\max\{|i|\,\big|\,a_{1,i}\cdot b_{1,i}\cdot d_{1,i}\cdot
e_{1,i}\ne0\},\, \ \,M_{1,2}=\max\{|r|\,\big|\,
b_{1,r}\cdot c_{1,r}\cdot f_{1,r}\cdot g_{1,r}\ne0\}.
\end{eqnarray*}
Using the induction on $M_{1,1}+M_{1,2}$,
and replacing $d_0$ by $d_0 - u_{\rm inn}$, where $u$ is a
combination of some $L_{i}\otimes L_{-i}$, $G_{p}\otimes G_{-p}$,
$T_{r}\otimes T_{-r}$, $G_{r}\otimes T_{-r}$, $T_{r}\otimes G_{-r}$,
$L_{j}\otimes G_{-j}$ and $G_{j}\otimes L_{-j}$, one can suppose
\begin{eqnarray*}
&&a_{1,i}=b_{1,p}=c_{1,r}=d_{1,j} =e_{1,k}=f_{1,r_1}=g_{1,r_2}=0,
\end{eqnarray*}
for any $i\in\Z\bs\{-2,1\}$, $j\in\Z\bs\{1\}$, $k\in\Z\bs\{-2\}$,
$p\in\frac{1}{2}\Z\bs\{-\frac{3}{2},\frac{1}{2}\}$,
$r\in\frac{1}{2}+\Z$, $r_1\in\frac{1}{2}+\Z\bs\{{\frac{1}{2}}\}$,
 $r_2\in\frac{1}{2}+\Z\bs\{-{\frac{3}{2}}\}$,
using which one can rewrite $d_0(L_{1})$ as
\begin{eqnarray}\label{0502n02}
d_0(L_{1})\!\!\!&=\!\!\!&a_{1,-2}L_{-1}\otimes
L_{2}+a_{1,1}L_{2}\otimes
L_{-1}+b_{1,-\frac{3}{2}}G_{-\frac{1}{2}}\otimes
G_{\frac{3}{2}}+b_{1,\frac{1}{2}}G_{\frac{3}{2}}\otimes
G_{-\frac{1}{2}}\nonumber\\
&&+d_{1,1}L_{2}\otimes G_{-1}+e_{1,-2}G_{-1}\otimes
L_{2}+f_{1,\frac{1}{2}}G_{\frac{3}{2}}\otimes
T_{-\frac{1}{2}}+g_{1,-{\frac{3}{2}}}T_{-{\frac{1}{2}}}\otimes
G_{{\frac{3}{2}}}.
\end{eqnarray}
Then one can see that $2L_{-1}*d_0(L_{1})$ is equal to
\begin{eqnarray*}
&&-6a_{1,-2}L_{-1}\otimes
L_{1}-6a_{1,1}L_{1}\otimes L_{-1}-6d_{1,1}L_{1}\otimes
G_{-1}-6e_{1,-2}G_{-1}\otimes
L_{1}\\
&&+d_{1,1}L_{2}\otimes
G_{-2}+e_{1,-2}G_{-2}\otimes
L_{2}-4b_{1,-\frac{3}{2}}G_{-\frac{1}{2}}\otimes
G_{\frac{1}{2}}-4b_{1,\frac{1}{2}}G_{\frac{1}{2}}\otimes
G_{-\frac{1}{2}}\\
&&-4f_{1,\frac{1}{2}}G_{\frac{1}{2}}\otimes
T_{-\frac{1}{2}}+f_{1,\frac{1}{2}}G_{\frac{3}{2}}\otimes
T_{-\frac{3}{2}}+g_{1,-{\frac{3}{2}}}T_{-{\frac{3}{2}}}\otimes
G_{{\frac{3}{2}}}-4g_{1,-{\frac{3}{2}}}T_{-{\frac{1}{2}}}\otimes
G_{{\frac{1}{2}}},
\end{eqnarray*}
while $L_1*d_0(L_{-1})-\mbox{$\sum\limits_{p\in\frac{1}{2}\Z}$}\big((\frac{3}{2}-p)b_{-1,p}
+(\frac{3}{2}+p)b_{-1,p+1}\big)G_{p}\otimes G_{-p}$ is equal to
\begin{eqnarray*}
&&\!\!\!\!\!\!\!\!\mbox{$\sum\limits_{i\in\Z}$}\big((2-i)a_{-1,i}+(2+i)a_{-1,i+1}\big)L_{i}\otimes
L_{-i}+\mbox{$\sum\limits_{r\in\frac{1}{2}+\Z}$}\!\big((1-r)c_{-1,r}+(1+r)c_{-1,r+1}\big)T_{r}\otimes
T_{-r}\\
&&\!\!\!\!\!\!\!\!+\mbox{$\sum\limits_{j\in\Z}$}\Big(\!\big((2-j)d_{-1,j}
+(\frac{3}{2}+j)d_{-1,j+1}\big)L_{j}\otimes
G_{-j}+\big((\frac{3}{2}-j)e_{-1,j}
+(2+j)e_{-1,j+1}\big)G_{j}\otimes L_{-j}\Big)
\\
&&\!\!\!\!\!\!\!\!+\mbox{$\sum\limits_{r\in\frac{1}{2}+\Z}$}\Big(\!\big(\!(\frac{3}{2}-r)f_{-1,r}
+(r+1)f_{-1,r+1})\!\big)G_{r}\otimes
T_{-r}+\big(\!(1-r)g_{-1,r}
+(\frac{3}{2}+r)g_{-1,r+1}\big)T_{r}\otimes G_{-r}\!\Big).
\end{eqnarray*}
Applying $d_0$ to $[L_{1},L_{-1}] = 2L_0$ and using Claim
\ref{sub1}, we obtain
\begin{eqnarray}\label{0516a01}
L_{-1}*d_0(L_{1})=L_1*d_0(L_{-1}).
\end{eqnarray}
Comparing the coefficients of $T_{r}\otimes T_{-r}$, $L_{i}\otimes
L_{-i}$ and $G_{p}\otimes G_{-p}$ in (\ref{0516a01}), one has
\begin{eqnarray*}
&&(r-1)c_{-1,r}
=(r+1)c_{-1,r+1},\\
&&3a_{-1,-1}+a_{-1,0}+3a_{1,-2} =a_{-1,1}+3a_{-1,2}+3a_{1,1}=
(i-2)a_{-1,i}-(i+2)a_{-1,i+1}=0,\\
&&2b_{-1,-\frac{1}{2}}+b_{-1,\frac{1}{2}}+2b_{1,-\frac{3}{2}}
=b_{-1,\frac{1}{2}}+2b_{-1,\frac{3}{2}}+2b_{1,\frac{1}{2}}=
(p-\frac{3}{2})b_{-1,p}-(p+\frac{3}{2})b_{-1,p+1}=0,
\end{eqnarray*}
for any $i\in\Z\bs\{\pm1\}$, $r\in\frac{1}{2}+\Z$ and
$p\in\frac{1}{2}\Z\bs\{\pm\frac{1}{2}\}$, which together with our
suppose that all the sets $\{a_{-1,i}\,|\,i\in\Z\}$,
$\{b_{-1,p}\,|\,p\in\frac{1}{2}\Z\}$,
$\{c_{-1,r}\,|\,r\in\frac{1}{2}+\Z\}$ are of finite rank, imply
\begin{eqnarray*}
&&c_{-1,r}=b_{-1,p}=2b_{-1,-\frac{1}{2}}+b_{-1,\frac{1}{2}}+2b_{1,-\frac{3}{2}}
=b_{-1,\frac{1}{2}}+2b_{-1,\frac{3}{2}}+2b_{1,\frac{1}{2}}=0,\\
&&a_{-1,i}=3a_{-1,-1}+a_{-1,0}+3a_{1,-2}
=a_{-1,1}+3a_{-1,2}+3a_{1,1}=a_{-1,0}+a_{-1,1}=0,
\end{eqnarray*}
for any $i\in\Z\bs\{\pm1,0,2\}$, $r\in\frac{1}{2}+\Z,$ and
$p\in\frac{1}{2}\Z\bs\{\pm\frac{1}{2},\frac{3}{2}\}$.
Comparing the coefficients of $L_{j}\otimes G_{-j}$ and
$G_{k}\otimes L_{-k}$ in (\ref{0516a01}), one has
\begin{eqnarray*}
&&(j-2)d_{-1,j}=(j+\frac{3}{2})d_{-1,j+1}, \ \
d_{-1,1}+\frac{5}{2}d_{-1,2}+3d_{1,1}
=\frac{7}{2}d_{-1,3}-\frac{1}{2}d_{1,1}=0,\\
&&(k-\frac{3}{2})e_{-1,k}=(k+2)e_{-1,k+1},\ \
\frac{5}{2}e_{-1,-1}+e_{-1,0}+3e_{1,-2}=
\frac{7}{2}e_{-1,-2}-\frac{1}{2}e_{1,-2}=0,
\end{eqnarray*}
for any $j\in\Z\bs\{1,2\}$ and $k\in\Z\bs\{-1,-2\}$, which together
with our suppose that all the sets $\{d_{-1,j}\,|\,j\in\Z\}$ and
$\{e_{-1,k}\,|\,k\in\Z\}$ are of finite rank, imply
\begin{eqnarray*}
&&d_{1,1}=e_{1,-2}=d_{-1,j}=e_{-1,j}=0,\ \ \forall\,\,j\in\Z.
\end{eqnarray*}
Comparing the coefficients of $G_r\otimes T_{-r}$ and
$T_{r}\otimes G_{-r}$ in (\ref{0516a01}), one has
\begin{eqnarray*}
&&(r-\frac{3}{2})f_{-1,r}=(1+r)f_{-1,r+1},\ \
f_{-1,\frac{1}{2}}+\frac{3}{2}f_{-1,\frac{3}{2}}+2f_{1,\frac{1}{2}}
=\frac{5}{2}f_{-1,\frac{5}{2}}-\frac{1}{2}f_{1,\frac{1}{2}}=0,\\
&&(s-1)g_{-1,s}=(\frac{3}{2}+s)g_{-1,s+1},\ \
\frac{3}{2}g_{-1,-\frac{1}{2}}+g_{-1,\frac{1}{2}}+2g_{1,-{\frac{3}{2}}}=
\frac{5}{2}g_{-1,-\frac{3}{2}}-\frac{1}{2}g_{1,-\frac{3}{2}}=0,
\end{eqnarray*}
for any $r\in\frac{1}{2}+\Z\bs\{\frac{1}{2},\frac{3}{2}\}$ and
$s\in\frac{1}{2}+\Z\bs\{-\frac{1}{2},-\frac{3}{2}\}$, which together
with our suppose that all the sets
$\{f_{-1,r}\,|\,r\in\frac{1}{2}+\Z\}$ and
$\{g_{-1,r}\,|\,r\in\frac{1}{2}+\Z\}$ are of finite rank, give
\begin{eqnarray*}
&&f_{1,\frac{1}{2}}=g_{1,-\frac{3}{2}}=0=f_{-1,r}=g_{-1,r},\ \
\forall\,\,r\in\frac{1}{2}+\Z.
\end{eqnarray*}
Then one can rewrite $d_0(L_{-1})$ $\big({\rm see}
(\ref{0502n01})\big)$ as
\begin{eqnarray*}
d_0(L_{-1})\!\!\!&=\!\!\!&3(a_{-1,2}\!+\!a_{1,1})(L_{-1}\!\otimes L_{0}\!-\!L_{0}\!\otimes
L_{-1})\!+\!a_{-1,2}L_{1}\!\otimes
L_{-2}\!-\!(a_{-1,2}+a_{1,1}+a_{1,-2})L_{-2}\!\otimes
L_{1}\\
&&+(b_{-1,\frac{3}{2}}+b_{1,\frac{1}{2}}-b_{1,-\frac{3}{2}})G_{-\frac{3}{2}}\!\otimes
G_{\frac{1}{2}}
-2(b_{-1,\frac{3}{2}}+b_{1,\frac{1}{2}})G_{-\frac{1}{2}}\otimes
G_{-\frac{1}{2}}+b_{-1,\frac{3}{2}}G_{\frac{1}{2}}\otimes
G_{-\frac{3}{2}}.
\end{eqnarray*}
Applying $d_0$ to $[L_{2},L_{-1}]=3L_1$, we obtain
\begin{eqnarray}\label{0517a01}
L_{2}*d_0(L_{-1})-L_{-1}*d_0(L_{2})=3d_0(L_1),
\end{eqnarray}
where
\begin{eqnarray*}
&&\!\!\!\!\!\!\!L_{2}\ast
d_0(L_{-1})=-4(a_{-1,2}+a_{1,1}+a_{1,-2})L_{0}\otimes
L_{1}-(a_{-1,2}+a_{1,1}+a_{1,-2})L_{-2}\otimes
L_{3}\\
&&\!\!\!\!\!\!\!+3(a_{-1,2}+a_{1,1})(3L_{1}\otimes
L_{0}-2L_{2}\otimes L_{-1})+3(a_{-1,2}+a_{1,1})(2L_{-1}\otimes
L_{2}-3L_{0}\otimes
L_{1})\\
&&\!\!\!\!\!\!\!+a_{-1,2}L_{3}\otimes L_{-2}+4a_{-1,2}L_{1}\otimes
L_{0}+(b_{-1,\frac{3}{2}}+b_{1,\frac{1}{2}}-b_{1,-\frac{3}{2}})(\frac{5}{2}G_{\frac{1}{2}}\otimes
G_{\frac{1}{2}}+\frac{1}{2}G_{-\frac{3}{2}}\otimes
G_{\frac{5}{2}})\\
&&\!\!\!\!\!\!\!
-3(b_{-1,\frac{3}{2}}+b_{1,\frac{1}{2}})(G_{\frac{3}{2}}\otimes
G_{-\frac{1}{2}}+G_{-\frac{1}{2}}\otimes G_{\frac{3}{2}})
+\frac{1}{2}b_{-1,\frac{3}{2}}G_{\frac{5}{2}}\otimes G_{-\frac{3}{2}}
+\frac{5}{2}b_{-1,\frac{3}{2}}G_{\frac{1}{2}}\otimes
G_{\frac{1}{2}},
\end{eqnarray*}
and $L_{-1}\ast d_0(L_2)$ can be rewritten as
\begin{eqnarray*}
&&\!\!\!\!\!\!
\mbox{$\sum\limits_{j}$}\big((j-\frac{3}{2})d_{2,j-1}-(j+3)d_{2,j}
\big)L_{j+1}\otimes G_{-j}+\mbox{$\sum\limits_{k}$}\big((k-2)e_{2,k-1}-(k+\frac{5}{2})e_{2,k}
\big)G_{k+1}\otimes L_{-k}\\
&&\!\!\!\!\!\!+\mbox{$\sum\limits_{r}$}\big((r-\frac{3}{2})g_{2,r-1}-(r+2)g_{2,r}
\big)T_{r+1}\otimes G_{-r}
+\mbox{$\sum\limits_{r}$}\big((r-1)f_{2,r-1}\!-\!(r+\frac{5}{2})f_{2,r}
\big)G_{r+1}\otimes T_{-r}\\
&&\!\!\!\!\!\!+\mbox{$\sum\limits_{i}$}
\big((i-2)a_{2,i-1}-(i+3)a_{2,i}\big)L_{i+1}\otimes L_{-i}
+\mbox{$\sum\limits_{p}$}\big((p-\frac{3}{2})b_{2,p-1}
-(p+\frac{5}{2})b_{2,p}\big)G_{p+ 1}\otimes G_{-p}\\
&&\!\!\!\!\!\!+\mbox{$\sum\limits_{r}$}\big((r-1)c_{2,r-1}
-(r+2)c_{2,r}\big)T_{r+1}\otimes T_{-r}.
\end{eqnarray*}
Comparing the coefficients of $G_{p+ 1}\otimes G_{-p}$ and
$T_{r+1}\otimes T_{-r}$ in (\ref{0517a01}) with
$p\in\Z,\,r\in\frac{1}{2}+\Z$, we obtain
\begin{eqnarray*}
&&(p-\frac{3}{2})b_{2,p-1}=(p+\frac{5}{2})b_{2,p},\ \
(r-1)c_{2,r-1}=(r+2)c_{2,r},
\end{eqnarray*}
which together with our suppose that all the sets
$\{b_{2,p}\,|\,p\in\Z\}$ and $\{c_{2,r}\,|\,r\in\frac{1}{2}+\Z\}$
are of finite rank, give
\begin{eqnarray*}
&&b_{2,p}=c_{2,r}=0,\ \ \forall\,\,p\in\Z,\,r\in\frac{1}{2}+\Z.
\end{eqnarray*}
Comparing the coefficients of $L_{i+1}\otimes L_{-i}$ in
(\ref{0517a01}) with $i\in\Z$, one has
\begin{eqnarray*}
&&\mbox{$\sum\limits_{i\in\Z}$}
\big((i-2)a_{2,i-1}-(i+3)a_{2,i}\big)L_{i+1}\otimes L_{-i}
+3(2a_{-1,2}+3a_{1,1})L_{2}\otimes L_{-1}\\
&&=-(a_{-1,2}+a_{1,1}+a_{1,-2})L_{-2}\otimes
L_{3}+3(2a_{-1,2}+2a_{1,1}-a_{1,-2})L_{-1}\otimes
L_{2}\\
&&-(13a_{-1,2}+13a_{1,1}+4a_{1,-2})L_{0}\otimes
L_{1}+(13a_{-1,2}+9a_{1,1})L_{1}\otimes L_{0}+a_{-1,2}L_{3}\otimes L_{-2},
\end{eqnarray*}
which together with our suppose that the set
$\{a_{2,i}\,|\,i\in\Z\}$ is of finite rank, give
\begin{eqnarray*}
&&a_{1,-2}+a_{1,1}=a_{-1,2}=a_{2,i}=0,\ \
4a_{2,-3} =-15a_{1,1}-a_{2,0},\\
&&a_{2,-2}=6a_{1,1}+a_{2,0},\ \ 2a_{2,-1}=-9a_{1,1}-3a_{2,0},\ \ \
4a_{2,1}=9a_{1,1}-a_{2,0},
\end{eqnarray*}
for any $i\in\Z\bs\{-3,\cdots,1\}$. Comparing the coefficients of $G_{p+1}\otimes G_{-p}$ in
(\ref{0517a01}) with $p\in\frac{1}{2}+\Z$, one has
\begin{eqnarray*}
&&\mbox{$\sum\limits_{p}$}\big((p-\frac{3}{2})b_{2,p-1}
-(p+\frac{5}{2})b_{2,p}\big)G_{p+1}\otimes G_{-p}
-\frac{1}{2}b_{-1,\frac{3}{2}}G_{\frac{5}{2}}\otimes
G_{-\frac{3}{2}}\\
&&=\frac{1}{2}(b_{-1,\frac{3}{2}}+b_{1,\frac{1}{2}}-b_{1,-\frac{3}{2}})G_{-\frac{3}{2}}\otimes
G_{\frac{5}{2}}-3(b_{-1,\frac{3}{2}}+b_{1,\frac{1}{2}}
+b_{1,-\frac{3}{2}})G_{-\frac{1}{2}}\otimes G_{\frac{3}{2}}\\
&&\ \ \
+\frac{5}{2}(2b_{-1,\frac{3}{2}}+b_{1,\frac{1}{2}}-b_{1,-\frac{3}{2}})G_{\frac{1}{2}}\otimes
G_{\frac{1}{2}}-3(b_{-1,\frac{3}{2}}+2b_{1,\frac{1}{2}})G_{\frac{3}{2}}\otimes
G_{-\frac{1}{2}},
\end{eqnarray*}
which together with our suppose that the set
$\{b_{2,p}\,|\,p\in\frac{1}{2}+\Z\}$ is of finite rank, give
\begin{eqnarray*}
&&b_{2,-\frac{3}{2}}=-6b_{1,\frac{1}{2}}+3b_{2,\frac{1}{2}}
=-b_{2,-\frac{1}{2}},\ \,b_{1,-\frac{3}{2}}=b_{1,\frac{1}{2}},\ \,b_{2,-\frac{5}{2}}=
4b_{1,\frac{1}{2}}-b_{2,\frac{1}{2}},\ \,b_{-1,\frac{3}{2}}=b_{2,p}=0,
\end{eqnarray*}
for any $p\in\frac{1}{2}+\Z\bs\{-\frac{5}{2},\cdots,\frac{1}{2}\}$. Comparing the coefficients of $L_{j+1}\otimes G_{-j}$ in
(\ref{0517a01}) with $j\in\Z$, one has
\begin{eqnarray*}
&&\mbox{$\sum\limits_{j\in\Z}$}\big((j-\frac{3}{2})d_{2,j-1}-(j+3)d_{2,j}
\big)L_{j+1}\otimes G_{-j}=0,
\end{eqnarray*}
which together with our suppose that the set
$\{d_{2,j}\,|\,j\in\Z\}$ is of finite rank, give
\begin{eqnarray*}
d_{2,j}=0,\ \forall\,\,j\in\Z.
\end{eqnarray*}
Comparing the coefficients of $G_{k+1}\otimes L_{-k}$ in
(\ref{0517a01}) with $k\in\Z$, one has
\begin{eqnarray*}
&&\mbox{$\sum\limits_{k\in\Z}$}\big((k-2)e_{2,k-1}-(k+\frac{5}{2})e_{2,k}
\big)G_{k+1}\otimes L_{-k}=0,
\end{eqnarray*}
which together with our suppose that the set
$\{e_{2,k}\,|\,k\in\Z\}$ is of finite rank, give
\begin{eqnarray*}
e_{2,k}=0,\ \forall\,\,k\in\Z.
\end{eqnarray*}
Comparing the coefficients of $G_{r+1}\otimes T_{-r}$ in
(\ref{0517a01}) with $r\in\frac{1}{2}+\Z$, one has
\begin{eqnarray*}
&&\mbox{$\sum\limits_{r\in\frac{1}{2}+\Z}$}\big((r-1)f_{2,r-1}\!-\!(r+\frac{5}{2})f_{2,r}
\big)G_{r+1}\otimes T_{-r}=0,
\end{eqnarray*}
which together with our suppose that the set
$\{f_{2,r}\,|\,r\in\frac{1}{2}+\Z\}$ is of finite rank, give
\begin{eqnarray*}
f_{2,r}=0,\ \forall\,\, r\in\frac{1}{2}+\Z.
\end{eqnarray*}
Comparing the coefficients of $T_{r+1}\otimes G_{-r}$ in
(\ref{0517a01}) with $r\in\frac{1}{2}+\Z$, one has
\begin{eqnarray*}
&&\mbox{$\sum\limits_{r\in\frac{1}{2}+\Z}$}\big((r-\frac{3}{2})g_{2,r-1}-(r+2)g_{2,r}
\big)T_{r+1}\otimes G_{-r}=0,
\end{eqnarray*}
which together with our suppose that the set
$\{f_{2,r}\,|\,r\in\frac{1}{2}+\Z\}$ is of finite rank, give
\begin{eqnarray*}
g_{2,r}=0,\ \ \forall\,\, r\in\frac{1}{2}+\Z.
\end{eqnarray*}
Then one can rewrite $d_0(L_{1})$, $d_0(L_{-1})$ and $d_0(L_{2})$ as
\begin{eqnarray*}
d_0(L_{1})\!\!\!&=&\!\!\!a_{1,1}(L_{2}\otimes
L_{-1}-L_{-1}\otimes L_{2})+b_{1,\frac{1}{2}}(G_{-\frac{1}{2}}\otimes
G_{\frac{3}{2}}+G_{\frac{3}{2}}\otimes
G_{-\frac{1}{2}}),\\[4pt]
d_0(L_{-1})\!\!\!&=&\!\!\!3a_{1,1}(L_{-1}\otimes L_{0}-L_{0}\otimes
L_{-1})-2b_{1,\frac{1}{2}}G_{-\frac{1}{2}}\otimes
G_{-\frac{1}{2}},\\[4pt]
d_0(L_{2})\!\!\!&=&\!\!\!-\frac{1}{4}(15a_{1,1}+a_{2,0})L_{-1}\otimes
L_{3}+(6a_{1,1}+a_{2,0})L_{0}\otimes
L_{2}-\frac{3}{2}(3a_{1,1}+a_{2,0})L_{1}\otimes
L_{1}\\
&&\!\!\!+a_{2,0}L_{2}\otimes
L_{0}+\frac{1}{4}(9a_{1,1}-a_{2,0})L_{3}\otimes
L_{-1}+(4b_{1,\frac{1}{2}}-b_{2,\frac{1}{2}})G_{-\frac{1}{2}}\otimes
G_{\frac{5}{2}}\\[4pt]
&&\!\!\!+b_{2,\frac{1}{2}}G_{\frac{5}{2}}\otimes
G_{-\frac{1}{2}}+3(2b_{1,\frac{1}{2}}-b_{2,\frac{1}{2}})(G_{\frac{3}{2}}\otimes
G_{\frac{1}{2}}-G_{\frac{1}{2}}\otimes G_{\frac{3}{2}}).
\end{eqnarray*}
Applying $d_0$ to $[L_1,L_{-2}]=3L_{-1}$, one has
\begin{eqnarray}\label{0627a01}
L_1\ast d_0(L_{-2})-L_{-2}\ast  d_0(L_1)=3d_0(L_{-1}),
\end{eqnarray}
while $L_1\ast d_0(L_{-2})$ can be written as
\begin{eqnarray*}
&&\!\!\!\!\!\!\!\!\mbox{$\sum\limits_{j}$}\big((3-j)d_{-2,j}
+(j+\frac{3}{2})d_{-2,j+1}\big)L_{j-1}\otimes G_{-j}
\!+\!\mbox{$\sum\limits_{k}$}\big((\frac{5}{2}-k)e_{-2,k}
+(k+2)e_{-2,k+1}\big)G_{k-1}\otimes
L_{-k}\\
&&\!\!\!\!\!\!\!\!+\mbox{$\sum\limits_{r}$}\big((\frac{5}{2}-r)f_{-2,r}
+(r+1)f_{-2,r+1}\big)G_{r-1}\otimes T_{-r}\!+\!\mbox{$\sum\limits_{r}$}\big((2-r)g_{-2,r}
+(r+\frac{3}{2})g_{-2,r+1}\big)T_{r-1}\otimes G_{-r},\\
&&\!\!\!\!\!\!\!\!+\mbox{$\sum\limits_{i}$}\big((3-i)a_{-2,i}
+(i+2)a_{-2,i+1}\big)L_{i-1}\otimes
L_{-i}\!+\!\mbox{$\sum\limits_{p}$}\big((\frac{5}{2}-p)b_{-2,p}
+(p+\frac{3}{2})b_{-2,p+1}\big)G_{p-1}\otimes G_{-p}\\
&&\!\!\!\!\!\!\!\!+\mbox{$\sum\limits_{r}$}\big((2-r)c_{-2,r}
+(r+1)c_{-2,r+1}\big)T_{r-1}\otimes
T_{-r},
\end{eqnarray*}
and
\begin{eqnarray*}
L_{-2}\ast  d_0(L_1)\!\!\!&=&\!\!\!a_{1,1}(L_{-3}\otimes
L_{2}-L_{2}\otimes L_{-3})+4a_{1,1}(L_{-1}\otimes
L_{0}-L_{0}\otimes L_{-1})\\
&&\!\!\!-\frac{1}{2}(b_{1,\frac{1}{2}}G_{-\frac{5}{2}}\otimes
G_{\frac{3}{2}}+b_{1,\frac{1}{2}}G_{\frac{3}{2}}\otimes
G_{-\frac{5}{2}})-5b_{1,\frac{1}{2}}G_{-\frac{1}{2}}\otimes
G_{-\frac{1}{2}}.
\end{eqnarray*}
Comparing the coefficients of $L_{i-1}\otimes L_{-i}$ in
(\ref{0627a01}), one has
\begin{eqnarray*}
&&\mbox{$\sum\limits_{i\in\Z}$}\big((3-i)a_{-2,i}
+(i+2)a_{-2,i+1}\big)L_{i-1}\otimes
L_{-i}\\
&&=a_{1,1}L_{-3}\otimes L_{2}+13a_{1,1}L_{-1}\otimes
L_{0}-13a_{1,1}L_{0}\otimes L_{-1}-a_{1,1}L_{2}\otimes L_{-3},
\end{eqnarray*}
which forces
\begin{eqnarray*}
&&a_{1,1}=a_{-2,i}=0,\ \ \forall\,\,i\in\Z.
\end{eqnarray*}
Comparing the coefficients of $G_{p-1}\otimes G_{-p}$ in
(\ref{0627a01}), one has
\begin{eqnarray*}
&&\mbox{$\sum\limits_{p}$}\big((\frac{5}{2}-p)b_{-2,p}
+(p+\frac{3}{2})b_{-2,p+1}\big)G_{p-1}\otimes G_{-p}\\
&&=-\frac{1}{2}b_{1,\frac{1}{2}}G_{-\frac{5}{2}}\otimes
G_{\frac{3}{2}}-11b_{1,\frac{1}{2}}G_{-\frac{1}{2}}\otimes
G_{-\frac{1}{2}}-\frac{1}{2}b_{1,\frac{1}{2}}G_{\frac{3}{2}}\otimes
G_{-\frac{5}{2}},
\end{eqnarray*}
which forces
\begin{eqnarray*}
&&b_{1,\frac{1}{2}}=b_{-2,p}=0,\ \ \forall\,\,p\in\frac{1}{2}\Z.
\end{eqnarray*}
Comparing the coefficients of $T_{r-1}\otimes T_{-r}$,
$L_{j-1}\otimes G_{-j}$, $G_{k-1}\otimes L_{-k}$, $G_{r-1}\otimes
T_{-r}$ and $T_{r-1}\otimes G_{-r}$ in (\ref{0627a01}), one has
\begin{eqnarray*}
&&\!\!\!\!\!\!(r-2)c_{-2,r}=(r+1)c_{-2,r+1},\,\
(j-3)d_{-2,j}=(j+\frac{3}{2})d_{-2,j+1},\,\
(k-\frac{5}{2})e_{-2,k}=(k+2)e_{-2,k+1},\\
&&\!\!\!\!\!\!(r-\frac{5}{2})f_{-2,r}=(r+1)f_{-2,r+1},\,\
(r-2)g_{-2,r}=(r+\frac{3}{2})g_{-2,r+1},
\end{eqnarray*}
which force
\begin{eqnarray*}
&&c_{-2,r}=d_{-2,j}=e_{-2,k}=f_{-2,r} =g_{-2,r}=0,\ \
\forall\,\,j,k\in\Z,\,r\in\frac{1}{2}+\Z.
\end{eqnarray*}
Then one can rewrite $d_0(L_{\pm1})$ and $d_0(L_{\pm2})$ as
\begin{eqnarray*}
d_0(L_{1})\!\!\!&=&\!\!\!d_0(L_{-1})=d_0(L_{-2})=0,\\[4pt]
d_0(L_{2})\!\!\!&=&\!\!\!-(\frac{1}{4}a_{2,0}L_{-1}\otimes
L_{3}+a_{2,0}L_{3}\otimes
L_{-1})+a_{2,0}(L_{0}\otimes L_{2}+L_{2}\otimes L_{0})
-\frac{3}{2}a_{2,0}L_{1}\otimes
L_{1}\\[4pt]
&&\!\!\!+b_{2,\frac{1}{2}}(G_{\frac{5}{2}}\otimes
G_{-\frac{1}{2}}-G_{-\frac{1}{2}}\otimes
G_{\frac{5}{2}})-3b_{2,\frac{1}{2}}(G_{\frac{3}{2}}\otimes
G_{\frac{1}{2}}-G_{\frac{1}{2}}\otimes G_{\frac{3}{2}}).
\end{eqnarray*}
Applying $d_0$ to $[L_{2},L_{-2}]=4L_{0}$, one has
$L_{-2}\ast  d_0(L_{2})=0$, which implies
\begin{eqnarray*}
&&a_{2,0}=b_{2,\frac{1}{2}}=0.
\end{eqnarray*}
Then one can claim $d_0(L_{2})$ is equal to zero. Hence
$d_0(L_{i})=0$, $\forall\,\,i\in\Z$.

One can write
\begin{eqnarray*}
d_0(G_{\frac{1}{2}})\!\!\!&=\!\!\!&\mbox{$\sum\limits_{i\in\Z}$}(a_{i}L_{i}\otimes
T_{\frac{1}{2}-i} +b_{i}T_{\frac{1}{2}-i}\otimes L_{i}
+c_{i}L_{i}\otimes G_{\frac{1}{2}-i}
+d_{i}G_{\frac{1}{2}-i}\otimes L_{i}\\
&&e_{i}G_{i}\otimes T_{\frac{1}{2}-i} +f_{i}T_{\frac{1}{2}-i}\otimes
G_{i}+g_{i}G_{i}\otimes G_{\frac{1}{2}-i}
+h_{i}G_{\frac{1}{2}-i}\otimes G_{i}),
\end{eqnarray*}
where the sums are all finite. Applying $d_0$ to
$[L_1,G_{\frac{1}{2}}]=0$ and comparing the coefficients, one can
deduce $d_0(G_{\frac{1}{2}})$ must be zero, which together with
$[L_m,G_{\frac{1}{2}}]=\frac{m-1}{2}G_{\frac{1}{2}+m}$ $(\forall\,\,
m\in\Z)$ and $[L_{-1},G_{\frac{5}{2}}]=-3G_{\frac{3}{2}}$ forces
\begin{eqnarray}\label{0812n01}
d_0(G_{r})=0,\ \ \forall\,\,r\in\frac{1}{2}+\Z.
\end{eqnarray}

Write $d_0(G_{0})$ as
\begin{eqnarray*}
&&\mbox{$\sum\limits_{i\in\Z}$}a_{i}L_{i}\otimes
L_{-i}+\mbox{$\sum\limits_{p\in\frac{1}{2}\Z}$}b_{p}G_{p}\otimes
G_{-p}+\mbox{$\sum\limits_{r\in\frac{1}{2}+\Z}$}c_{r}T_{r}\otimes
T_{-r}+\mbox{$\sum\limits_{j\in\Z}$}d_{j}L_{j}\otimes
G_{-j}\nonumber\\
&&+\mbox{$\sum\limits_{k\in\Z}$}e_{k}G_{k}\otimes
L_{-k}+\mbox{$\sum\limits_{r\in\frac{1}{2}+\Z}$}f_{r}G_{r}\otimes
T_{-r}+\mbox{$\sum\limits_{r\in\frac{1}{2}+\Z}$}g_{r}T_{r}\otimes
G_{-r},
\end{eqnarray*}
where the sums are all finite. Applying $d_0$ to $[G_0,G_0]=2L_0$
and comparing the coefficients, one can deduce
(\,$\forall\,\,r\in\frac{1}{2}+\Z,\ i\in\Z$)
\begin{eqnarray*}
&&rb_{r}=c_{r},\ \ f_{r}=g_{r},\ \ ia_{i}=-4b_{i},\ \ d_{i}=-e_{i}.
\end{eqnarray*}
Then $d_0(G_0)$ can be rewritten as
\begin{eqnarray*}
d_0(G_{0})\!\!\!&=\!\!\!&\mbox{$\sum\limits_{i\in\Z}$}(a_{i}L_{i}\otimes
L_{-i}-\frac{ia_{i}}{4}G_{i}\otimes G_{-i}+d_{i}L_{i}\otimes
G_{-i}-d_{i}G_{i}\otimes L_{-i})\\
&&+\mbox{$\sum\limits_{r\in1/2+\Z}$}(b_{r}G_{r}\otimes
G_{-r}+rb_{r}T_{r}\otimes T_{-r}+f_{r}G_{r}\otimes
T_{-r}+f_{r}T_{r}\otimes G_{-r}).
\end{eqnarray*}
Applying $d_0$ to $[L_{-1},[L_{1},G_{0}]]=-\frac{3}{4}G_0$, we
obtain
\begin{eqnarray}\label{0807m01}
&&L_{-1}\ast L_{1}\ast d_0(G_0)=-\frac{3}{4}d_0(G_0).
\end{eqnarray}
Comparing the coefficients of $L_{i}\otimes L_{-i}$, $L_{i}\otimes
G_{-i}$, $G_{r}\otimes G_{-r}$, $G_{r}\otimes T_{-r}$ and
$T_{r}\otimes T_{-r}$, one has
\begin{eqnarray*}
&&(8i^2-13)a_{i}=4(i-2)^2a_{i-1}+4(i+2)^2a_{i+1},\\
&&(4i^2-4)d_{i}=(i-2)(2i-3)d_{i-1}+(i+2)(2i+3)d_{i+1},\\
&&2r^2f_{r}+(r-1)f_{r-1}(\frac{3}{2}-r)
+(r+1)(-r-\frac{3}{2})f_{r+1}=0,\\
&&(r-\frac{3}{4})b_{r}+(\frac{3}{2}-r)(r-\frac{3}{2})b_{r-1}
+(\frac{3}{2}+r)(-r-\frac{3}{2})b_{r+1}=0,
\end{eqnarray*}
which together with our suppose the set $\{i\,|\,i\in\Z,\ a_i\cdot
b_{i+\frac{1}{2}}\cdot d_i\neq0\}$ is of finite rank, imply
\begin{eqnarray*}
&&a_{i}=0,\ \ \ b_{r}=f_{r}=0,\ \ \ d_{j}=0,\ \ \ d_{-1}=-d_1,
\end{eqnarray*}
for any $i\in\Z$, $r\in\frac{1}{2}+\Z$ and $j\in\Z\bs\{\pm1\}$. Then
$d_0(G_0)$ can be rewritten as
\begin{eqnarray}\label{0812n02}
d_0(G_{0})=-d_{1}L_{-1}\otimes G_{1}+d_{1}L_{1}\otimes
G_{-1}+d_{1}G_{-1}\otimes L_{1}-d_{1}G_{1}\otimes L_{-1}.
\end{eqnarray}

One can write $d_0(T_{\frac{1}{2}})$ as
\begin{eqnarray*}
d_0(T_{\frac{1}{2}})\!\!\!&=\!\!\!&\mbox{$\sum\limits_{i\in\Z}$}a_{i}L_{i}\otimes
T_{\frac{1}{2}-i}
+\mbox{$\sum\limits_{i\in\Z}$}b_{i}T_{\frac{1}{2}-i}\otimes L_{i}
+\mbox{$\sum\limits_{i\in\Z}$}c_{i}L_{i}\otimes G_{\frac{1}{2}-i}
+\mbox{$\sum\limits_{i\in\Z}$}d_{i}G_{\frac{1}{2}-i}\otimes L_{i}\\
&&+\mbox{$\sum\limits_{i\in\Z}$}e_{i}G_{i}\otimes T_{\frac{1}{2}-i}
+\mbox{$\sum\limits_{i\in\Z}$}f_{i}T_{\frac{1}{2}-i}\otimes G_{i}
+\mbox{$\sum\limits_{i\in\Z}$}g_{i}G_{i}\otimes G_{\frac{1}{2}-i}
+\mbox{$\sum\limits_{i\in\Z}$}h_{i}G_{\frac{1}{2}-i}\otimes G_{i},
\end{eqnarray*}
where the sums are all finite. Applying $d_0$ to
$[L_{-1},[L_{1},T_{\frac{1}{2}}]]=\frac{3}{4}T_{\frac{1}{2}}$, we
obtain
\begin{eqnarray}\label{0812a1}
&&L_{-1}\ast L_{1}\ast
d_0(T_{\frac{1}{2}})=\frac{3}{4}d_0(T_{\frac{1}{2}}).
\end{eqnarray}
Comparing the coefficients of $L_{i}\otimes T_{\frac{1}{2}-i}$ and
$T_{\frac{1}{2}-i}\otimes L_{i}$ in the both sides of
(\ref{0812a1}), we obtain
\begin{eqnarray*}
&&(2i^2-i-2)x_{i}+(2-i)(i-\frac{3}{2})x_{i-1}
-(i+\frac{1}{2})(i+2)x_{i+1}=0\ \ \mbox{for}\ \ x=a\ \mbox{or}\ b,
\end{eqnarray*}
which imply
\begin{eqnarray*}
&&x_{i}=0,\ \ \forall\,\,i\in\Z\bs\{0,\pm1\}\ \ \mbox{while}\ \
x_0=-2x_1=-2x_{-1}\ \ \mbox{for}\ \ x=a\ \mbox{or}\ b.
\end{eqnarray*}
Comparing the coefficients of $L_{i}\otimes G_{\frac{1}{2}-i}$ and
$G_{\frac{1}{2}-i}\otimes L_{i}$ in the both sides of
(\ref{0812a1}), we obtain
\begin{eqnarray*}
&&(2i^2-i-\frac{11}{4})y_{i}-(i-2)^2y_{i-1}-(i+1)(i+2)y_{i+1}=0\ \
\mbox{for}\ \ y=c\ \mbox{or}\ d,
\end{eqnarray*}
which imply
\begin{eqnarray*}
&&y_{i}=0\ \ \mbox{for}\ \ y=c\ \mbox{or}\ d,\ \forall\,\,i\in\Z.
\end{eqnarray*}
Comparing the coefficients of $G_{i}\otimes T_{\frac{1}{2}-i}$ and
$T_{\frac{1}{2}-i}\otimes G_{i}$ in the both sides of
(\ref{0812a1}), we obtain
\begin{eqnarray*}
&&(2i^2-i-\frac{3}{4})z_{i}-(\frac{3}{2}-i)^2z_{i-1}-(i+\frac{1}{2})(i+\frac{3}{2})z_{i+1}=0\
\ \mbox{for}\ \ z=e\ \mbox{or}\ f,
\end{eqnarray*}
which imply
\begin{eqnarray*}
&&z_{i}=0\ \ \mbox{for}\ \ z=e\ \mbox{or}\ f,\ \forall\,\,i\in\Z.
\end{eqnarray*}
Comparing the coefficients of $G_{i}\otimes G_{\frac{1}{2}-i}$ and
$G_{\frac{1}{2}-i}\otimes G_{i}$ in the both sides of
(\ref{0812a1}), we obtain
\begin{eqnarray*}
&&(2i^2-i-\frac{3}{2})w_{i}+(\frac{3}{2}-i)(i-2)w_{i-1}-(i+1)(i+\frac{3}{2})w_{i+1}=0\
\ \mbox{for}\ \ w=g\ \mbox{or}\ h,
\end{eqnarray*}
which imply
\begin{eqnarray*}
&&w_{i}=0,\ \ \forall\,\,i\in\Z\bs\{0,1\}\ \ \mbox{while}\ \
w_0=-w_1\ \ \mbox{for}\ \ w=g\ \mbox{or}\ h.
\end{eqnarray*}
Then $d_0(T_{\frac{1}{2}})$ can be rewritten as
\begin{eqnarray}
d_0(T_{\frac{1}{2}})\!\!\!&=\!\!\!&a_{1}(L_{-1}\otimes
T_{\frac{3}{2}}-2L_{0}\otimes T_{\frac{1}{2}}+L_{1}\otimes
T_{-\frac{1}{2}})+g_{1}(G_{1}\otimes G_{-\frac{1}{2}}-G_{0}\otimes
G_{\frac{1}{2}})+\nonumber\\
&&b_{1}(T_{\frac{3}{2}}\otimes L_{-1}-2T_{\frac{1}{2}}\otimes
L_{0}+T_{-\frac{1}{2}}\otimes L_{1})+h_{1}(G_{-\frac{1}{2}}\otimes
G_{1}-G_{\frac{1}{2}}\otimes G_{0}).\label{0812n03}
\end{eqnarray}
Applying $d_0$ to
$[G_0,G_{\frac{1}{2}}]=\frac{1}{2}T_{\frac{1}{2}}$, one has
\begin{eqnarray*}
&&2G_0\cdot d_0(G_{\frac{1}{2}})+2G_{\frac{1}{2}}\cdot
d_0(G_{0})=d_0(T_{\frac{1}{2}}).
\end{eqnarray*}
Using (\ref{0812n01}), (\ref{0812n02}) and (\ref{0812n03}) and
comparing the coefficients of all the products in the above
identity, we obtain
\begin{eqnarray*}
&&d_{1}=a_{1}=g_{1}=b_{1}=h_{1}=0,
\end{eqnarray*}
which implies
\begin{eqnarray*}
&&d_0(G_0)=d_0(T_{\frac{1}{2}})=0.
\end{eqnarray*}
By now, we have proved
\begin{eqnarray*}
&&d_0(L_i)=d_0(G_r)=d_0(G_0)=d_0(T_{\frac{1}{2}})=0,\ \
\forall\,\,i\in\Z,\ r\in\frac{1}{2}+\Z,
\end{eqnarray*}
which together with (\ref{LieB}), implies $d_0(\LL)=0$ for the case
$d_0\in\mbox{Der}_{\bar 0}(\LL,\VV)$. Thus the claim follows.

\begin{clai}\adddot\label{sub2+}
Suppose $d_0\in\Der_{\bar1}(\WW,\VV)$ is odd. By
replacing $d_0$ by $d_0-u_{\rm inn}$ for some $u\in \VV_0$, we can
suppose $d_0(\LL)= 0 $.
\end{clai}
Employing the similar techniques used in Claim \ref{sub2}, one can
see the claim holds.
\begin{clai}\adddot
The sum in (\ref{summable}) is finite.
\end{clai}
For any $i\in\Z$, suppose $d_{i}=(v_{i})_{\rm inn}$ for some
$v_{i}\in \VV_{i}$. If $|\{i\,|\,v_i\ne0\}|$ is infinite, then
$d(L_0)=\sum_{i\in\Z}L_0\ast v_{i}=-\sum_{i\in\Z}i v_{i}$ is an
infinite sum, which contradicts $d\in\mbox{Der}(\LL,\VV)$. Thus the
claim and proposition follow. \QED
\begin{lemm}\adddot \label{lemma3.4} If $r\in \VV$ satisfies
$x\ast  r\in {\rm Im}(\bo\otimes \bo-\tau)\,(\,\forall\,\,x\in
\LL)$, then $r\in {\rm Im}(\bo\otimes \bo-\tau)$.
\end{lemm}
\ni{\it Proof.~}~Note $\LL \ast  {\rm Im}(\bo\otimes\bo-\tau)\subset
{\rm Im}(\bo\otimes\bo-\tau).$ Write $r=\sum_{i\in\frac{1}{2}\Z}r_i$
with $r_i\in \VV_i$. Obviously, $r\in{\rm Im}(\bo\otimes\bo-\tau)$
if and only if $r_i\in{\rm Im}(\bo\otimes\bo-\tau)$ for all
$i\in\frac{1}{2}\Z.$ Thus without loss of generality, one can
suppose $r=r_i$ is homogeneous.

If $i\in\frac{1}{2}\Z^*$, then $r_i=-\frac1i L_0\ast  r_i\in{\rm
Im}(\bo\otimes\bo-\tau)$. For the case $i=0$, one can write
\begin{eqnarray*}
r_0\!\!\!&=\!\!\!&\mbox{$\sum\limits_{i\in\Z}$}a_{i}L_{i}\otimes
L_{-i}+\mbox{$\sum\limits_{p\in\frac{1}{2}\Z}$}b_{p}G_{p}\otimes
G_{-p}+\mbox{$\sum\limits_{r\in\frac{1}{2}+\Z}$}c_{r}T_{r}\otimes
T_{-r}+\mbox{$\sum\limits_{j\in\Z}$}d_{j}L_{j}\otimes
G_{-j}\nonumber\\
&&+\mbox{$\sum\limits_{k\in\Z}$}e_{k}G_{k}\otimes
L_{-k}+\mbox{$\sum\limits_{r\in\frac{1}{2}+\Z}$}f_{r}G_{r}\otimes
T_{-r}+\mbox{$\sum\limits_{r\in\frac{1}{2}+\Z}$}g_{r}T_{r}\otimes
G_{-r},
\end{eqnarray*}
where the sum are all finite.  Since the elements of the form
$u_{1,i}:=L_i\otimes L_{-i}-L_{-i}\otimes L_{i}$,
$u_{2,p}:=G_p\otimes G_{-p}-G_{-p}\otimes G_{p}$,
$u_{3,r}:=T_r\otimes T_{-r}-T_{-r}\otimes T_{r}$,
$v_{i}:=L_i\otimes G_{-i}-G_{-i}\otimes L_{i}$ and
$w_{r}:=G_r\otimes T_{-r}-T_{-r}\otimes G_{r}$ are all in ${\rm
Im}(\bo\otimes \bo-\tau),$ replacing $v$ by $v-u$, where $u$ is a combination of
some $u_{1,i}$, $u_{2,p}$, $u_{3,r}$, $v_{i}$  and $w_{r}$, one can suppose
\begin{eqnarray}
&&a_{i}\ne 0\ \,\Longrightarrow\ \,i\in\Z_+,\label{wpqr2}\\
&&b_{p}\ne 0\ \,\Longrightarrow\ \,p\in\frac{1}{2}\Z_+,\label{wpqr3}\\
&&c_{r}\ne 0\ \,\Longrightarrow\ \,r\in\frac{1}{2}+\Z_+,\label{wpqr4}\\
&&e_i=g_r=0,\ \ \forall\,\,i\in\Z,\ r\in\frac{1}{2}+\Z.\label{wpqr1}
\end{eqnarray}
Then $r_0$ can be rewritten as
\begin{eqnarray*}
r_0\!\!\!\!&=\!\!\!\!&\mbox{$\sum\limits_{i\in\Z_+}$}\!\!a_{i}L_{i}\otimes
L_{-i}+\!\!\mbox{$\sum\limits_{p\in\frac{1}{2}\Z_+}$}\!\!b_{p}G_{p}\otimes
G_{-p}+\!\!\mbox{$\sum\limits_{r\in\frac{1}{2}+\Z_+}$}\!\!c_{r}T_{r}\otimes
T_{-r}+\!\!\mbox{$\sum\limits_{j\in\Z}$}\!\!d_{j}L_{j}\otimes
G_{-j}+\!\!\mbox{$\sum\limits_{r\in\frac{1}{2}+\Z}$}\!\!f_{r}G_{r}\otimes
T_{-r}.
\end{eqnarray*}
First assume that $a_{i}\ne 0$ for some $i>0$. Choose $j<0$ such
that $i+j<0$. Then we see that the term $L_{i}\otimes L_{j-i}$
appears in $L_{j}\cdot r_0,$ but (\ref{wpqr2}) implies that the term
$L_{i+j}\otimes L_{-i}$ does not appear in $L_j\cdot r_0$, a
contradiction with the fact that $L_j\cdot r_0\in {\rm Im}(\bo\otimes \bo-\tau)$.
Then one further can suppose $a_i=0,\ \forall\,\,i\in\Z^*$.
Similarly, one also can suppose $b_p=c_r=0$ for all $p\in\frac{1}{2}\Z^*$, $r\in\frac{1}{2}+\Z$.
Therefore, $r_0$ can be rewritten as
\begin{eqnarray}\label{sm2}
r_0=\mbox{$\sum\limits_{j\in\Z}$}d_{j}L_{j}\otimes
G_{-j}+\mbox{$\sum\limits_{r\in\frac{1}{2}+\Z}$}f_{r}G_{r}\otimes
T_{-r}+a_{0}L_0\otimes L_{0}+b_{0}G_0\otimes G_{0}\,.
\end{eqnarray}
Finally, we mainly use the fact ${\rm Im}(\bo\otimes \bo-\tau)\subset{\rm
Ker}(\bo\otimes \bo+\tau)$ and the assumption that $\LL\cdot r_0\subset{\rm
Im}(\bo\otimes \bo-\tau)$ to deduce $a_{0}=d_{0}=d_j=f_r=0$ for all $j\in\Z$, $r\in\frac{1}{2}+\Z$. One
has
\begin{eqnarray*}
0\!\!\!&=&\!\!\!(\bo\otimes \bo+\tau)(L_1\cdot r_0)\\
\!\!\!&=&\!\!\!2a_{0}(L_1\otimes L_{0}+L_1\otimes L_{0})+b_{0}(G_1\otimes G_{0}+G_0\otimes G_{1})\\
&&\!\!\!+\mbox{$\sum\limits_{r\in1/2+\Z}$}
\big((3/2-r)f_{r-1}+rf_{r}\big)(G_{r}\otimes T_{1-r}+T_{1-r}\otimes G_{r})\\
&&\!\!\!+\mbox{$\sum\limits_{j\in\Z}$}\big((2-j)d_{j-1}+(1/2+j)d_{j}\big)
(L_{j}\otimes G_{1-j}+G_{1-j}\otimes L_{j})\,.
\end{eqnarray*}
Then noticing both the sets $\{j\,|\,d_j\ne0\}$ and $\{r\,|\,f_r\ne0\}$ of finite rank and
comparing the coefficients of the tensor products, one immediately
gets
\begin{eqnarray*}
&&a_0=b_0=d_j=f_r=0,\ \
\forall\,\,p\in\Z,\ r\in\frac{1}{2}+\Z.
\end{eqnarray*}
Thus the lemma follows.\QED
\vskip10pt

\ni{\it Proof of Theorem \ref{mainthe}.} Let $(\LL
,[\cdot,\cdot],\D)$ be a Lie super-bialgebra structure on $\LL$. Then $\D=\D_r$ is
defined by (\ref{e-D-r}) for some $r\in\VV_{\bar0} $.
By (\ref{cond2}), ${\rm Im}\, \D\subset{\rm
Im}(\bo\otimes\bo-\tau)$. Thus by Lemma \ref{lemma3.4}, $r\in{\rm
Im}(1\otimes 1-\tau)$. Then  (\ref{cond2}), (\ref{e-p2.1}) and Corollary
\ref{colo} show that $c(r)=0$. Thus $(\LL,[\cdot,\cdot],\D)$ is triangular coboundary. \QED

\end{document}